\documentclass[12pt,reqno,a4paper]{article}

\usepackage{amsmath,amsthm,amstext,amscd,amssymb,euscript,mathrsfs,enumitem}
\usepackage{calrsfs}
\usepackage{epsfig}
\usepackage{hyperref,tikz,relsize}

\newcommand{\Z}{\mathbb Z}
\newcommand{\R}{\mathbb R}
\newcommand{\N}{\mathbb N}

\newcommand{\E}{\mathbb E}
\newcommand{\Zd}{\mathbb Z^d}

\renewcommand{\phi}{\varphi}

\newcommand{\pee}{\ensuremath{\mathbb{P}}}

\def\1{{\mathchoice {\rm 1\mskip-4mu l} {\rm 1\mskip-4mu l}
{\rm 1\mskip-4.5mu l} {\rm 1\mskip-5mu l}}}

\newtheorem{theorem}{{\small T}{\scriptsize HEOREM}}[section]
\newtheorem{corollary}{{\bf{\small C}{\scriptsize OROLLARY}}}[section]
\newtheorem{proposition}{{\bf{\small P}{\scriptsize ROPOSITION}}}[section]
\newtheorem{lemma}{{\bf{\small L}{\scriptsize EMMA}}}[section]
\newtheorem{remark}{{\bf{\small R}{\scriptsize EMARK}}}[section]
\newtheorem{definition}{{\bf{\small D}{\scriptsize EFINITION}}}[section]

\renewenvironment{proof}[1]
{\noindent{{\bf{\small{ P}{\scriptsize ROOF}}}.}\hspace{0.1cm} #1} {$\;\qed$\newline}

\newcommand{\beq}{\begin{eqnarray}}
\newcommand{\eeq}{\end{eqnarray}}

\newcommand{\ba}{\begin{align*}}
\newcommand{\ea}{\end{align*}}

\newcommand{\be}{\begin{equation}}
\newcommand{\ee}{\end{equation}}

\newcommand{\bl}{\begin{lemma}}
\newcommand{\el}{\end{lemma}}

\newcommand{\br}{\begin{remark}}
\newcommand{\er}{\end{remark}}

\newcommand{\bt}{\begin{theorem}}
\newcommand{\et}{\end{theorem}}

\newcommand{\bd}{\begin{definition}}
\newcommand{\ed}{\end{definition}}

\newcommand{\bp}{\begin{proposition}}
\newcommand{\ep}{\end{proposition}}

\newcommand{\bc}{\begin{corollary}}
\newcommand{\ec}{\end{corollary}}

\newcommand{\bpr}{\begin{proof}}
\newcommand{\epr}{\end{proof}}

\newcommand{\bi}{\begin{itemize}}
\newcommand{\ei}{\end{itemize}}

\newcommand{\ben}{\begin{enumerate}}
\newcommand{\een}{\end{enumerate}}

\newcommand{\caD}{{\EuScript D}}

\newcommand{\caI}{{\mathcal I}}

\newcommand{\caN}{{\mathcal N}}

\newcommand{\caP}{{\mathcal P}}

\newcommand{\caS}{{\mathcal S}}

\usepackage{color}

\begin{document}
\title{Ergodic theory of the symmetric inclusion process}
\author{Kevin Kuoch$^{\textup{{\tiny 1 }}}$,
Frank Redig$^{\textup{{\tiny 2 }}}$
}

\footnotetext[1]{Johann Bernoulli Institute, Rijksuniversiteit Groningen. Postbus 407, 9700AK Groningen, The Netherlands. \\ E-mail: \texttt{kevin.kuoch@gmail.com}}
\footnotetext[2]{Delft Institute of Applied Mathematics, TU Delft. Mekelweg 4, 2628CD Delft, The Netherlands. \\ E-mail: \texttt{f.h.j.redig@tudelft.nl}}

\maketitle

\begin{abstract}
We prove the existence of a successful coupling  for $n$ particles in the symmetric inclusion process.
As a consequence we characterize the ergodic measures with finite moments, and obtain
sufficient conditions for a measure to converge  in the course of time to an invariant product measure.

\bigskip

\noindent

\end{abstract}

\section{Introduction}
In \cite[Chapter VIII]{Lig85}, a rather complete ergodic theory is given for the symmetric exclusion process (SEP).
In particular, for the simple symmetric exclusion process the only extremal invariant measures are
Bernoulli measures with constant density.
This complete characterization of the set of invariant measures is quite exceptional and in the case of SEP is a consequence of the fact the SEP is self-dual. Because of this, invariant measures
can be related to {\em bounded harmonic functions for the finite SEP}.
Then, by the construction of a successful coupling of the SEP with a finite number of particles, it is shown
that all bounded harmonic functions are constant, i.e. only depending on the number of particles.
From this in turn,  one can conclude that all invariant measures for the SEP are permutation invariant, from which
one derives by the De Finetti theorem that they are convex combinations of Bernoulli measures.
In \cite{GKR,GRV} an attractive version (in the sense of having attractive interaction between the particles)
of the SEP is introduced and called the {\em simple inclusion process} (SIP). In the SIP, particles perform nearest-neighbor jumps
according to a simple symmetric random walk and interact by ``inclusion jumps'', where pairs of neighboring particles jump to the same site at rate 1.
This analogy between SIP (attractive) and SEP (repulsive) becomes even more apparent in \cite{GRV}, where it is shown that the SIP satisfies the analogue of Liggett's comparison inequality \cite[Chapter VIII, Proposition I.7]{Lig85} for
the evolution of positive definite symmetric functions. The expectation at time $t>0$ of such a function
in the course of the evolution of $n$ SIP-particles is {\em larger than
in the course of the evolution of $n$ independent random walkers}. In particular this implies that a certain class of product measures is mapped by the evolution under the SIP to measures with
positive correlations (as opposed to negative correlations in the SEP).

In this paper, we want to go as far as possible in the study of the invariant measures of the SIP, i.e.,
understand its ergodic measures and their attractors. Because the number of particles is unbounded,
and because we want to use self-duality, we will have to restrict to a set of measures with all moments finite.
The main problem is then to construct a successful coupling for two sets of $n$ SIP-particles initially at different locations. That this coupling
seems possible is due to  the fact that as long as SIP-particles do not collide, i.e., are not at neighboring positions, they behave
as independent random walkers and these can be coupled by the coordinate-wise Ornstein coupling in any dimension.
The idea of the coupling of SIP-particles comes from \cite{dMP}, combined with \cite{OR}. In \cite{OR} it is shown
that inclusion particles and independent random  walkers can be coupled in such a way that at time
$t$ they are $o(\sqrt{t})$ apart. The period of time $[0,(1-\delta)t]$ in which coupling according to \cite{OR} is used (stage 1) is then followed
a period of time $[(1-\delta)t, t]$ (stage 2) in which the coordinate-wise Ornstein coupling
of independent random walkers is used {\em both for the independent walkers as well as for the SIP walkers}. The only problem then in order for this
coupling to be successful is to estimate the probability of being
coupled before a collision takes place. One can understand however that such a collision event is highly improbable (as $t\to\infty$),
because after a long time the walkers are much further apart ($O(\sqrt{t})$) than the distance between
the walkers and their inclusion partners ($o(\sqrt{t})$).
Once one has the successful coupling of SIP-particles, and as a consequence results on the structure of the invariant measures of the SIP,
all these results can be transferred without effort to corresponding results for interacting diffusion processes of which the SIP is
a dual process such as the BEP (Brownian Energy Process) and BMP (Brownian Momentum Process).

The rest of our paper is organized as follows. In Section \ref{2} we give basic definitions and
set up notations, in Section \ref{3} we prove the successful coupling, in Section \ref{4} we characterize
the class of ergodic so-called tempered (with finite moments) measures and in Section \ref{5}  we
give sufficient conditions  for a measure to converge to an invariant product measure.
\section{Notations and definitions}\label{2}
\subsection{The symmetric inclusion process}
We denote by $p(.,.)$ the transition probability of a simple symmetric nearest neighbor random walk
on the lattice $\Zd$, i.e.,
\[
p(x,y) =
\begin{cases}
\frac1{2d} \ \mbox{if}\ |x-y|=1\\
0\ \mbox{otherwise}. 
\end{cases}
\]
The {\em simple symmetric inclusion process} with parameter $m>0$ (denoted SIP($m$)) is an interacting particle
system where particles perform independent random walks according to the transition probabilities
$p(.,.)$ and on top of that, they interact by inclusion, i.e. each particle
``invites'' any other particle at nearest neighbor position at rate 1 to join its site (invitations are always followed up).
These ``invitation jumps'', or ``\textsl{inclusion jumps}'', create an attractive interaction between the particles.
This has to be compared with the interaction between particles of the symmetric exclusion process (SEP) where jumps joining two particles at the same site are {\em forbidden}. Here, on the contrary, these jumps are {\em encouraged}.

More formally,
the SIP is a continuous-time Markov process $(\eta_t)_{t \geq 0} \in \mathbb N^S$ whose generator $\mathcal L$ acts on local functions $f$, i.e. depending only on a finite number of occupation variables, as
\begin{equation}
\mathcal L f (\eta) = \sum\limits_{x \in \Zd} \sum\limits_{y: \|y-x\|=1} p(x,y)\eta(x) \Big( \frac{m}{2} + \eta(y)\Big) \big( f(\eta^{x,y}) - f(\eta) \big),
\end{equation}
where $\eta^{x,y}\in \N^{\Zd}$ stands for the configuration obtained from $\eta\in \N^{\Zd}$ by moving one particle from $x$ to $y$, i.e. $\eta^{x,y} = \eta - \delta_x +\delta_y$, where $\delta_x$ denotes the configuration with a single particle at $x$ and no particles elsewhere ; $\| \cdot \|$ stands for the $\ell_1$-norm.

\subsection{Invariant product measures}
For $\lambda \in [0,1)$, define the homogeneous discrete Gamma product measure $\nu_\lambda^m$ on $\N^{\Zd}$ whose marginals are given by
\begin{equation}\label{def:refmeasure}
\nu_\lambda^m \big\{ \eta(x)=k \big\} = \dfrac{1}{Z_{\lambda,m}} \dfrac{\lambda^k}{k!} \dfrac{\Gamma (\frac{m}{2}+k)}{\Gamma(\frac{m}{2})}
\end{equation}
where $\Gamma(\cdot)$ denotes the Gamma function and
\[
Z_{\lambda,m} =\left(\frac{1}{1-\lambda}\right)^{m/2}
\]
is the normalizing constant. These measures $\nu_\lambda^m$ are reversible and ergodic for the SIP($m$),
see \cite{GRV} for more details. One of the main questions answered in the present paper is whether these are the only ergodic measures
within a certain class.
\subsection{Duality}
A duality relation is a link between a dual process with the process of interest in such a way it allows to perform computations for one process in terms of another.  The link is created via the duality function.
For further details see \cite[Chapter II, Section 4]{Lig85}.

\begin{definition}[duality relation]
Suppose $(\xi_t)_{t \geq 0}$ and $(\eta_t)_{t \geq 0}$ are Markov processes on $S_1$ and $S_2$ respectively. Let $D$ be a bounded measurable function on $S_1 \times S_2$. The processes $(\xi_t)_{t \geq 0}$ and $(\eta_t)_{t \geq 0}$ are said to be {\em dual} to one another with respect to $D$ if
\begin{equation}\label{bibi}
\mathbb E_\xi D(\xi_t,\eta) = \mathbb E_\eta D(\xi,\eta_t).
\end{equation}
\end{definition}
If the processes $(\xi_t)_{t \geq 0}$ and $(\eta_t)_{t \geq 0}$ are the same, then we call
\eqref{bibi} {\em self-duality.}
In that sense, the SIP($m$) is self-dual (see \cite {GKR}) with duality functions $D(\cdot,\cdot)$ given by
\begin{equation}\label{def:dualf}
D(\xi,\eta) = \prod\limits_{x \in S} d(\xi(x),\eta(x))
\end{equation}
where
\[
d(k,l) = \left\{
\begin{array}{cl}
\dfrac{l !}{(l-k)!} \dfrac{\Gamma(\frac{m}{2})}{\Gamma (\frac{m}{2}+k)} & \mbox{for}\ k \leq l
\\
0  &  \mbox{for}\ k>l
\end{array}
\right.,
\]
and $d(0,0)=1$.
Self-duality of the SIP($m$) then means
\begin{equation}\label{okkie}
\mathbb E^{SIP(m)}_\eta D(\xi,\eta_t) = \mathbb E^{SIP(m)}_\xi D(\xi_t,\eta),
\end{equation}
where $\mathbb E^{SIP(m)}_\eta$ denotes the expectation of a SIP($m$) starting from an initial configuration $\eta$
and where $\xi$ is a finite configuration (i.e., having a finite number of particles).

The self-duality functions $D(.,.)$ and the reference measure $\nu_\lambda^m$ are naturally connected via
\begin{equation}\label{def:relation-dualmeasure}
\int D(\xi,\eta) d\nu_\lambda^m(\eta) = \Big(\dfrac{\lambda}{1-\lambda}\Big)^{|\xi|},
\end{equation}
where $|\xi|$ denotes the number of particles in the finite configuration $\xi$.
We refer the reader to \cite{GKR,GRV} for the proof of the self-duality \eqref{okkie} and further details and properties of the SIP. By duality relations, we derive as well related results for interacting diffusions that are  dual to the SIP: the Brownian momentum process and the Brownian energy process, see \cite{GRV} and references therein.

The main advantage of self-duality is that we can study the SIP with infinitely many particles
by studying the SIP with a finite number of particles. Indeed, to know the time dependent expectations
of the polynomials $D(\xi, \eta)$ it suffices to follow the evolution of the particles
in the finite configuration $\xi$, and the initial configuration $\xi$.

While exploiting the self-duality property, we will necessarily restrict to starting measures with finite moments. Let us denote by $\caP$ the set of all probability measures on the configuration space $\N^{\mathbb Z^d}$. We then consider the class of so-called {\em tempered probability measures} defined as follows:
\bd
\be\label{temp}
\caP_t=\left\{ \mu: \mu\in \caP: \forall n\in \N: \sup_{|\xi|=n} \int D(\xi,\eta) \mu(d\eta)=: c_n<\infty\right\}.
\ee
where $c_n$ satisfies the Carleman moment condition
\[
\sum_{n=1}^\infty c_n^{-1/n}=\infty
\]
ensuring the moments $\int D(\xi,\eta) \mu(d\eta)$ characterize uniquely the measure $\mu$ \cite{Kleiber}.
\ed
First, remark that by self-duality, a tempered measure remains tempered in the course of the evolution of
the SIP. Indeed, if $\mu\in \caP_t$ then, by conservation of the number of particles in the finite
SIP, denoting
\[
\sup_{|\xi|=n} \int D(\xi,\eta) \mu(d\eta)=c_n,
\]
we have, for $|\xi|=n$ and $t>0$,
\[
\int \E_\eta D(\xi,\eta_t) \mu(d\eta) =\E_\xi \int D(\xi_t,\eta) \mu(d\eta)\leq c_n,
\]
hence the time-evolved measure $\mu_t$ is tempered for all $t>0$, with the same dominating constants $c_n, n\in \N$ as the one of the starting
measure $\mu$.

We are interested in characterizing of the invariant measures which are ergodic for the SIP, and belong to
$\caP_t$. We call $\caI$ the set of invariant probability measures for the SIP, and
by $\caI_t$ the set $\caI\cap\caP_t$ of tempered invariant measures.
Furthermore, we call $\caI_e$ the set of extreme points of $\caI$, i.e.,
the set of invariant and ergodic probability measures for the SIP.

For a measure $\mu\in \caP_t$ we denote its $D$-transform by
\be\label{dtrans}
\hat{\mu}(\xi)=\int D(\xi,\eta) \mu(d\eta).
\ee
Here $\xi$ varies in the set of finite configurations, which we denote further by $\Omega_n$ if $\xi$ is a configuration which contains $n$ particles.
Note that $\xi\in\Omega_n$ such that $n=|\xi|$ can be identified with an $n$-tuple
$x_1,\ldots,x_n$ via
\be\label{iden}
\xi=\sum_{i=1}^n \delta_{x_i}.
\ee
As a consequence $\hat{\mu}$ can also be viewed as a symmetric function on $\cup_{n\in \N}(\Zd)^n$.

The following result is then a straightforward consequence
of the self-duality of  the SIP.
\bp\label{prop:harmo}
 A probability measure
$\mu\in \caI_t$ if and only if its $D$-transform $\hat{\mu}$ is bounded harmonic for the SIP, i.e.,
if and only if for all $t>0$
\be\label{harmo}
\E_\xi \hat{\mu} (\xi_t)=\hat{\mu}(\xi)
\ee
\ep
As a consequence, the study of ergodic measures in $\caI_t$ is reduced to
the problem of identifying the set of bounded harmonic functions for the SIP.
This is done via the construction of a successful coupling, which implies that bounded harmonic functions are constant.
\subsection{Bounded harmonic functions and successful coupling}
Via the identification
\eqref{iden}
we can see the evolution of $n$ SIP-particles
initially at positions $\mathbf x=(x_1,\ldots, x_n)$ as a process
$ \mathbf X^S (t)=(X_1^S(t), \ldots, X_n^S(t))$ on $(\Zd)^n$ so that for each $1 \leq i \leq n$, $X_i^S(t)$ keeps track of the location of the SIP-particle $i$ started from $x_i$. We denote by $\pee^{SIP}_{\mathbf x}$ its path
space measure.
A coupling of two copies of the SIP starting initially at different locations $\mathbf x=(x_1,\ldots, x_n),\mathbf y=(y_1,\ldots,y_n)$ is then defined
as  usual as a process $\{\mathbf X^S(t),\mathbf Y^S(t)), t\geq 0\}$ with $\mathbf X^S(0)= \mathbf x, \mathbf Y^S(0)= \mathbf y$ with first (resp. second)
marginals $\{\mathbf X^S(t), t\geq 0\}$, the SIP starting from $\mathbf x$ (resp. $\{\mathbf Y^S(t), t\geq 0\}$, the SIP starting from $\mathbf y$).
Its path space measure is then denoted by $\widehat \pee_{\mathbf x,\mathbf y}^{SIP}$ where the hat stands for the joint distribution in the coupling.
The coupling time is defined via
\begin{equation}
\tau= \inf \{t\geq 0: \mathbf X^S(s)=\mathbf Y^S(s)\ \forall s\geq t \}
\end{equation}
where by convention $\inf(\emptyset)=\infty$.
The coupling is successful if $\tau<\infty$ $\widehat \pee_{\mathbf x,\mathbf y}^{SIP}$-almost surely for all $(\mathbf x,\mathbf y)\in (\Z^d)^{2n}$.
It is well known (see e.g. \cite[Chapter 2]{Lig85}) the existence of a successful coupling implies that all bounded harmonic functions are constant.
\subsection{Diffusion processes related to the SIP}
The SIP is related via duality to the BEP, a system of interacting diffusions with state space $[0,\infty)^{\Zd}$ and to
the BMP, a system of interacting diffusions with state space $\R^{\Zd}$. See \cite{GKR} for more details.
This implies that many results on the invariant measures and characterization of ergodic measures can  be transferred to these processes.
This ``transference'' is a consequence of the following proposition.
\bp\label{transprop}
Assume $\{\zeta_t: t\geq 0\}$ is a Feller process on the state space $K^{\Zd}$ with $K$ a Polish space.
Assume that the process is dual to the SIP with duality function $\caD(\xi, \zeta)$.
Then let
\[
\caP^\caD_t= \left\{ \mu: \mu\in \caP: \forall n\in \N: \sup_{|\xi|=n} \int_{K^{\Zd}}\caD(\xi,\zeta) \mu(d\zeta)<\infty\right\}
\]
denote the corresponding set of tempered probability measures, and
denote for $\mu\in \caP^\caD_t$, $\hat{\mu}(\xi)=\int \caD(\xi, \zeta)\mu(d\zeta)$ its $\caD$-transform.
Then we have $\mu\in \caP^\caD_t$ is invariant for $\{\zeta_t:t\geq 0\}$ if and only if
$\hat{\mu}$ is a bounded harmonic function for the SIP.
\ep
\section{Successful coupling for the SIP}\label{3}
In this section we prove:
\bt\label{success}
There exists a successful coupling for the SIP.
\et

We introduce here and onwards the following notations for the evolution of different sets of particles: sets of $n$ independent random walkers, denoted IRW-particles, $\mathbf X^I(t)=(X_1^I(t),\ldots,X_n^I(t))$ and $\mathbf Y^I(t)=(Y_1^I(t),\ldots,Y_n^I(t))$ ; sets of $n$ SIP-particles $\mathbf X^S(t)=(X_1^S(t),\ldots,X_n^S(t))$ and $\mathbf Y^S(t)=(Y_1^S(t),\ldots,Y_n^S(t))$. We denote by $\mathbf x = (x_1, \ldots,x_n)$ (resp. $\mathbf y = (y_1, \ldots,y_n)$) the initial locations of $\mathbf X^I$ and $\mathbf X^S$ (resp. $\mathbf Y^I$ and $\mathbf Y^S$).

Prior to the proof, we first define the notion of collision of the process $(\mathbf X^S(\cdot),\mathbf Y^S(\cdot))$ at time $t>0$.
We say that {\em a collision} happens at time $t>0$ if
two SIP-particles belonging to {\em a same set} are at nearest-neighbor positions at time $t$, i.e.,
$$\{\exists i\neq j, 1 \leq i, j\leq n : | X_i^S(t) - X_j^S(t) | = 1\ \mbox{or}\ | Y_i^S(t) - Y_j^S(t) | = 1\} $$
for a set of SIP-particles $\mathbf X^S=(X_1^S,...,X_n^S)$.
Notice that particles $X^S_i(t)$  and $Y^S_j(t)$ at neighboring positions is {\em not considered as a collision}, i.e., collisions
only happen within the same set of particles.

\bpr
The proof of a successful coupling is twofold, we first consider the case when $d \geq 3$ and then, the case when $d \leq 2$.
\bigskip

{\bf Transient case: $d\geq 3$}.\\

We start with the simplest case $d\geq 3$  where the random walk $X(t)$ based on $p(.,.)$ is transient. More precisely we  have that
\[
\pee_{x}  (|X(t)|>1, \ \forall t\geq 0) =: H(  x)>0
\]
and $H(x)\to 1$ when $x\to\infty$. As a consequence, by the union of events bound,
with positive probability, $n$ IRW-particles $\mathbf X^I(t) = (X_1^I(t),...,X_n^I(t))$ starting from initial positions $\mathbf x = (x_1,\ldots,x_n)$ for which all $|x_i-x_j|>R$ ($i\not= j$) are large enough, will never collide.

If during a lapse of time no collision happens, then the IRW-particles and their corresponding coupled SIP-particles perform exactly the same jumps. It is only when IRW-particles collide that their corresponding SIP partners can behave differently.

 Assume all initial positions satisfying $|x_i-x_j|, |y_i-y_j|>R$ for all $1 \leq i \neq j \leq n$. Then with positive probability $p(R)$ the IRW-particles starting at $(x_1,\ldots,x_n)$ will never collide
and neither will the IRW-particles starting from $(y_1,\ldots,y_n)$. Moreover $p(R)\to 1$ when $R\to\infty$.
The two sets of IRW-particles can be coupled by the coordinate-wise Ornstein coupling.
Now we couple each set of IRW-particles with a corresponding set of SIP-particles via the coupling described in \cite[Theorem 3.2]{OR}, i.e.,
both sets perform the same random walk jumps, and inclusion jumps are only performed by the SIP-particles.

In this coupling, when no collision happens for both sets of IRW-particles, the corresponding two sets of SIP-particles, respectively starting from $(x_1,\ldots,x_n)$ and $(y_1,\ldots,y_n)$ as well, behave exactly as the IRW-particles, since IRW- and SIP-particles only behave differently when they collide and hence can be coupled.

Therefore, the two sets of SIP-particles
can be coupled with positive probability. To show that they can be coupled almost surely from any initial locations $\mathbf x$ and $\mathbf y$, let them first move for some time $T>0$ and
then start the coupling just described.
With probability $\pi(T)$, close to one (as $T\to\infty$) any two SIP- or IRW- particles will be at distance larger than  $\alpha(T)$ where $\alpha(T)\to\infty$ when $T\to\infty$, and
correspondingly, from then on they can be coupled with probability $p(\alpha(T))$.
Therefore, the probability that they cannot be coupled is bounded from above by $(1-\pi(T))+ (1-p(\alpha(T)))$ which tends to zero
as  $T\to\infty$.
\bigskip

{\bf Recurrent case: $d=1,2$}\\
To tackle the case $d=1$, we first give an outline. Note that in the two-dimensional case the same arguments hold and we therefore omit the proof for $d=2$.

We follow the line of thought of  \cite{dMP}.
The coupling proceeds in two stages:
\ben
\item First stage: in the time interval $[0,(1-\delta)t]$ (where $0<\delta<1$ is fixed), the two sets of IRW-particles make the same jumps and the SIP partners follow according to the coupling
of \cite{OR}. After this first stage, with probability close to one (as $t\to\infty$) any two different IRW-particle within the same set, as well as their corresponding SIP-particles partners will be at distance
$O(\sqrt{t})$. The distance between the independent random walkers and their corresponding SIP partner, on the contrary will be of order
$o(\sqrt{t})$.
\item Second stage: in the time interval $[(1-\delta)t, t]$ the two sets of IRW-particles as well as the two sets of SIP walkers
are coupled via the coordinate-wise Ornstein coupling. 
As long as the SIP particles do not collide, this is indeed a coupling, because then the
SIP particles behave as independent random walkers. If such a collision does happen, then we say that we have a failed coupling attempt.
\een
\textit{Stage 1}.
During the time interval $[0, (1-\delta)t]$,  couple the two sets of IRW-particles so that they perform the same jumps, thus, for any $s \leq (1-\delta)t$,
\begin{equation}\label{dist_II}
\sum\limits_{i=1}^n \big|(X_i^I-Y_i^I)(s)\big| = \sum\limits_{i=1}^n  | x_i-y_i | =: k_n.
\end{equation}

\begin{figure}[h!]
\centering
\begin{tikzpicture}[scale=0.8]
\draw [->,line width=1pt] (-.5,-.5) -- (-.5,-7.5);
\node at (-.5,0) {$\mathsmaller{0}$};
\node at (-.5,-8) {$\mathsmaller{(1-\delta) t}$};
\draw[help lines,line width=1pt,step=2,draw=black] (1,0) to (7,0);
\draw (2,0) node[circle,fill,inner sep=1pt,label=above:$\mathsmaller{x_1}$](){};
\draw (6,0) node[circle,fill,inner sep=1pt,label=above:$\mathsmaller{x_2}$](){};
\draw [red, xshift=0cm] plot [smooth, tension=1] coordinates { (6,0) (7,-3) (6,-4) (4.5,-5) (5.5,-6) (5,-8)};
\draw [red] plot [smooth, tension=1] coordinates { (2,0) (1,-3) (3,-4) (1,-5) (3,-6) (1,-8)};
\draw [blue] plot [smooth, tension=1] coordinates { (6,0) (7,-3) (6,-4) (4.5,-5) (5.5,-6) (5,-8)};
\draw [blue] plot [smooth, tension=1] coordinates { (2,0) (1,-3) (3,-4) (1,-5) (3,-6) (1,-8)};
\draw (1,-8) node[circle,fill,inner sep=1pt,label=above:$\mathsmaller{X_1^I(t-\delta t)}$](){};
\draw (5,-8) node[circle,fill,inner sep=1pt,label=above:$\mathsmaller{X_2^I(t-\delta t)}$](){};
\draw (1,-8) node[circle,fill,inner sep=1pt,label=below:$\mathsmaller{X_1^S(t-\delta t)}$](){};
\draw (5,-8) node[circle,fill,inner sep=1pt,label=below:$\mathsmaller{X_2^S(t-\delta t)}$](){};
\draw[help lines,line width=1pt,step=2,draw=black] (1,-8) to (7,-8);
\draw[help lines,line width=1pt,step=2,dashed] (8,.5) to (8,-8.5);
\draw[help lines,line width=1pt,step=2,draw=black] (9,0) to (15,0);

\draw (12,0) node[circle,fill,inner sep=1pt,label=above:$\mathsmaller{y_1}$](){};
\draw (14,0) node[circle,fill,inner sep=1pt,label=above:$\mathsmaller{y_2}$](){};
\draw [red] plot [smooth, tension=1] coordinates { (12,0) (11,-3) (13,-4) (11,-5) (13,-6) (11,-8)};
\draw [red] plot [smooth, tension=1] coordinates { (14,0) (15,-3) (13.5,-4) (12,-5) (14,-6) (13,-8)};
\draw [blue] plot [smooth, tension=1] coordinates { (12,0) (11,-3) (13,-4) (11,-5) (13,-6) (11,-8)};
\draw [blue] plot [smooth, tension=1] coordinates { (14,0) (15,-3) (13.5,-4) (12,-5) (13,-6) (12,-8)};

\draw (11,-8) node[circle,fill,inner sep=1pt,label=above:$\mathsmaller{Y_1^I(t-\delta t)}$](){};
\draw (11,-8) node[circle,fill,inner sep=1pt,label=below:](){};
\draw (10.2,-8) node[label=below:$\mathsmaller{Y_1^I(t-\delta t)}$](){};
\draw (13,-8) node[circle,fill,inner sep=1pt,label=above:$\mathsmaller{Y_2^I(t-\delta t)}$](){};
\draw (12,-8) node[circle,fill,inner sep=1pt,label=below:$\mathsmaller{Y_2^S(t-\delta t)}$](){};
\draw[help lines,line width=1pt,step=2,draw=black] (9,-8) to (15,-8);
\draw (1,0) node[circle,fill,inner sep=1pt,label=above:$\mathsmaller{A}$](){};
\draw (7,0) node[circle,fill,inner sep=1pt,label=above:$\mathsmaller{B}$](){};
\draw (9,0) node[circle,fill,inner sep=1pt,label=above:$\mathsmaller{A}$](){};
\draw (15,0) node[circle,fill,inner sep=1pt,label=above:$\mathsmaller{B}$](){};
\end{tikzpicture}
\caption{Stage 1. SIP-particles $\mathbf X^S$ (resp. $\mathbf Y^S$) in blue and IRW-particles $\mathbf X^I$ (resp. $\mathbf Y^I$) in red both start from $\mathbf x$ (resp. $\mathbf y$) on the same segment $[A,B]$ of $\mathbb Z$. Between each set of particles, the moves of the IRW-particles are coupled so that the distance in-between is constant to $|x_1-x_2|$ (resp. $|y_1-y_2|)$. Each SIP-particle is distant from an IRW-particle by $o(\sqrt t)$ being coupled thanks to the coupling given by \cite[Theorem 3.2]{OR}. If no collision occurs, IRW- and SIP- particles follow the same path (in purple). At time $(1-\delta) t$, the distance between any pair of IRW-particle is $O(\sqrt t )$.}
\label{fig1}
\end{figure}
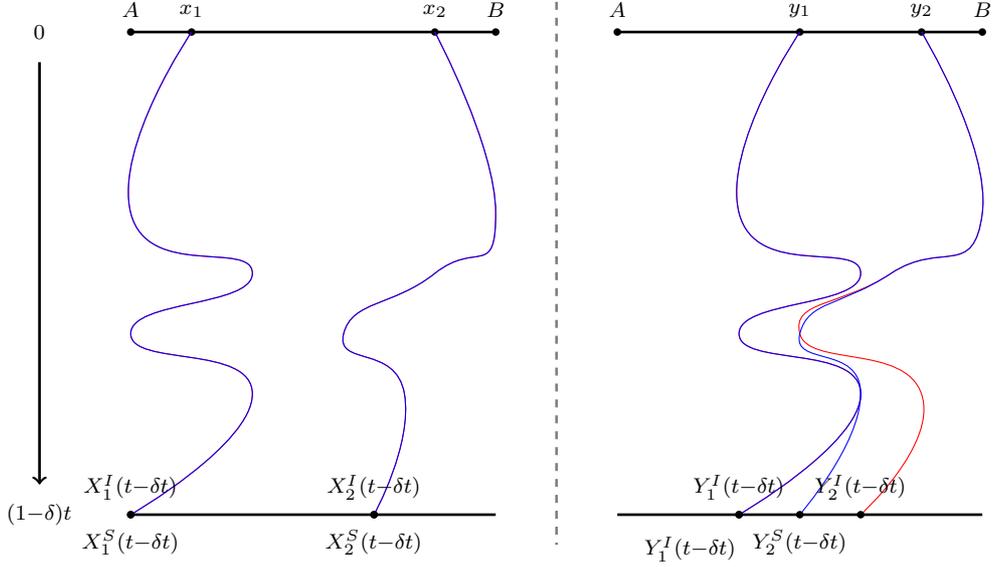

According to the coupling of \cite[Theorem 3.2]{OR}, two sets of SIP-particles $\mathbf X^S$ and $\mathbf Y^S$, starting respectively from $\mathbf x$ and $\mathbf y$, are coupled to two sets of IRW-particles $\mathbf X^I$ and $\mathbf Y^I$, starting respectively from $\mathbf x$ and $\mathbf y$ as well, such that their positions satisfy
\be\label{dist_IS}
\sum_{i=1}^n|X^S_i(s)-X^I_i(s)|+ |Y^S_i(s)-Y^I_i(s)|\leq \psi(s)
\ee
with
\be\label{mop}
\lim_{s\to\infty}\frac{\psi(s)}{\sqrt{s}}=0
\ee
with probability one (w.r.t. coupling distribution). By way of illustration, see Figure \ref{fig1}.

Gathering \eqref{dist_IS} and \eqref{dist_II}, it is now straightforward to see that for any $\delta \in (0,1)$, with high probability,
\begin{align*}
&\sum\limits_{i=1}^n  \big| \big( X_i^S - Y_i^S\big) (t - \delta t ) \big| \\
&\leq \sum\limits_{i=1}^n  \big| \big( X_i^S - X_i^I\big) (t - \delta t ) \big| +\big| \big( X_i^I - Y_i^I\big) (t - \delta t ) \big| + \big| \big( Y_i^I - Y_i^S\big) (t - \delta t ) \big|\\
& \leq k_n + 2 \psi(t-\delta t ).
\end{align*}

\textit{Stage 2}. Now, couple the two sets of IRW-particles and the SIP particles coordinate-wise during the time interval $[(1-\delta) t , t ]$ using the Ornstein coupling. If the two sets of SIP particles are coupled in this lapse of time $[t-\delta t , t ]$, and no collision occurred
then we say that the coupling attempt is succesful.

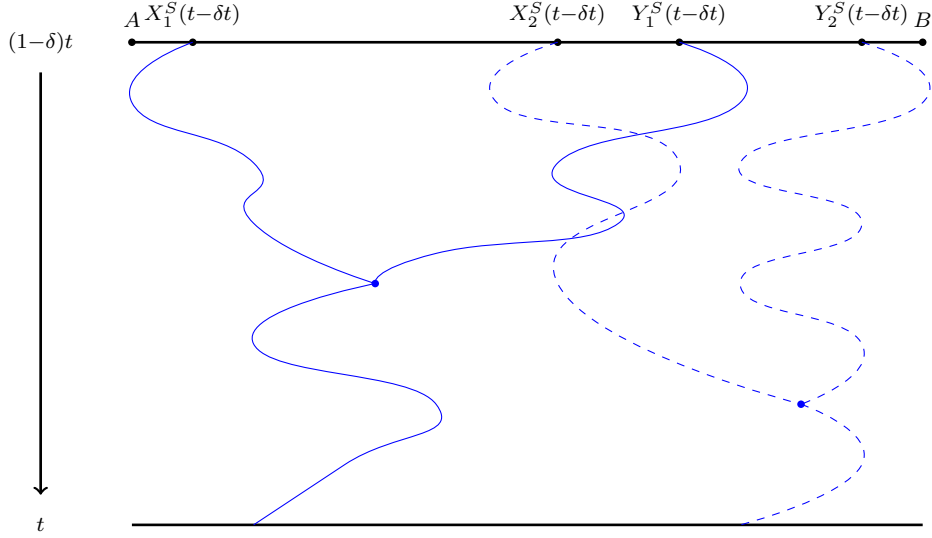
\begin{figure}[h!]
\centering
\begin{tikzpicture}[scale=0.8]
\draw [->,line width=1pt] (-.5,-.5) -- (-.5,-7.5);
\node at (-.5,0) {$\mathsmaller{(1-\delta) t}$};
\node at (-.5,-8) {$\mathsmaller{t}$};
\draw[help lines,line width=1pt,step=2,draw=black] (1,0) to (14,0);
\draw[help lines,line width=1pt,step=2,draw=black] (1,-8) to (14,-8);
\draw (1,0) node[circle,fill,inner sep=1pt,label=above:$\mathsmaller{A}$](){};
\draw (14,0) node[circle,fill,inner sep=1pt,label=above:$\mathsmaller{B}$](){};
\draw (10,0) node[circle,fill,inner sep=1pt,label=above:$\mathsmaller{Y_1^S(t-\delta t)}$](){};
\draw (2,0) node[circle,fill,inner sep=1pt,label=above:$\mathsmaller{X_1^S(t-\delta t)}$](){};
\draw (13,0) node[circle,fill,inner sep=1pt,label=above:$\mathsmaller{Y_2^S(t-\delta t)}$](){};
\draw (8,0) node[circle,fill,inner sep=1pt,label=above:$\mathsmaller{X_2^S(t-\delta t)}$](){};

\draw (5,-4) node[circle,fill, blue,inner sep=1pt](){};
\draw (12,-6) node[circle,fill, blue,inner sep=1pt](){};

\draw [blue] plot [smooth, tension=1] coordinates { (2,0) (1,-1) (3,-2) (3,-3) (5,-4)};
\draw [blue] plot [smooth, tension=1] coordinates { (10,0) (11,-1)  (8,-2) (9,-3) (6,-3.5) (5,-4)};
\draw [blue] plot [smooth, tension=1] coordinates { (5,-4) (3,-5) (6,-6) (4.5,-7) (3,-8)};
\draw [blue,dashed] plot [smooth, tension=1] coordinates { (8,0) (7,-1) (10,-2) (8,-4) (12,-6)};
\draw [blue,dashed] plot [smooth, tension=1] coordinates { (13,0) (14,-1) (11,-2) (13,-3) (11,-4) (13,-5) (12,-6)};
\draw [blue,dashed] plot [smooth, tension=1] coordinates { (12,-6) (13,-7) (11,-8)};

\end{tikzpicture}
\caption{Stage 2. While the IRW-particles $\mathbf X^I$ and $\mathbf Y^I$ are coupled coordinate-wise via an Ornstein coupling (there, omitted from the figure), the SIP-particles are consequently coupled in the same way provided no collision occurs within any set of particles (here, between $X_1^S$ and $X_2^S$ or $Y_1^S$ and $Y_2^S$, but any other does not matter).}
\label{fig2}
\end{figure}

After time $(1-\delta)t$, any pair of different SIP-particles as well as every pair of different IRW-particles of the same set (i.e., the ones starting from $\mathbf x$ as well as the ones starting from $\mathbf y$ are at distance of order $\sqrt{t}$ with probability close to one as $t\to\infty$. Indeed, for the IRW-particles this is clear from the invariance principle, whereas for the corresponding SIP-particles, it then follows via \eqref{dist_IS}. Therefore, it is sufficient to prove that for all $i$, the probability that $X^I(t)$ is coupled to $Y^I(t)$ before any collision
tends to zero as $t\to\infty$.  See Figure \ref{fig2}.

The motion of the difference of $X_i^I(t)-Y_i^I(t)$ in the Ornstein coupling is that of a simple continuous-time random walk at twice the speed (i.e., at rate $m= 2(m/2)$). So it suffices to see that for two IRW-particles moving at rate $m$, we have
\begin{equation}
\lim\limits_{t \rightarrow \infty} \mathbb {\widehat  P}_{\psi(t),\sqrt t}^{IRW} \big( \tau_0^1 \geq \tau_{\{-1,1\}}^2 \big) = 0,
\end{equation}
where $\tau_0^1$ (resp. $\tau^2_{\{-1,1\}}$) stands for the hitting time of $0$ (resp. of $\{-1,1\}$), $\mathbb {\widehat  P}_{a,b}^{IRW}$ denotes the joint distribution of two IRW-particles starting respectively from locations $a$ and $b$.

This is, using standard arguments, in turn implied by
\be\label{proba_collision}
\lim\limits_{t \rightarrow \infty} \mathbb {\widehat  P}_{\psi(t),\sqrt t}^{IRW} \big( \tau_0^1 \geq \tau_0^2 \big)=0
\ee
Recall by the reflection principle (see e.g. \cite{LL} Chapter I), for any $a \geq 0$:
\[
\mathbb P_0^{IRW} (\tau_a \geq t ) = \mathbb P^{IRW}_0 (|X_(t)| \leq a ),
\]
where here and in what follows $X(t)$ denotes a continuous-time simple random walk moving at rate $m$.
Therefore, choosing a time scale $\phi(t)$ such that $\phi(t)\to\infty$ as $t\to\infty$, we obtain
\begin{align*}
& \mathbb {\widehat P}_{\psi(t),\sqrt t}^{IRW} \Big( \tau_0^1 \geq \tau_0^2 \Big) \leq \mathbb {\widehat P}_{\psi(t),\sqrt t}^{IRW} \Big( \tau_0^1 \geq \tau_0^2 , \tau_0^2 \geq \phi(t)\Big)  + \mathbb {\widehat P}_{\psi(t),\sqrt t}^{IRW} \Big( \tau_0^1 \geq \tau_0^2 , \tau_0^2 \leq \phi(t)\Big)   \\
& \quad \leq \mathbb { P}_{\psi(t)}^{IRW} \Big( \tau_0^1 \geq \phi(t)\Big)  + \mathbb { P}_{\sqrt t}^{IRW} \Big(  \tau_0^2 \leq \phi(t)\Big)   \\
& \quad \leq \mathbb P^{IRW}_0 \Big( | X(\phi(t))| \leq \psi(t) \Big) + \mathbb P^{IRW}_0 \Big( | X(\phi(t)) | \geq \sqrt t \Big) \\
 & \quad = \mathbb P^{IRW}_0 \Bigg( \dfrac{| X({\phi(t)})|}{ \psi(t)} \leq 1 \Bigg) + \mathbb P^{IRW}_0 \Bigg( \dfrac{| X({\phi(t)}) | }{\sqrt t} \geq 1 \Bigg).
\end{align*}
Since $| X(\phi(t))| \sim \phi(t)^{1/2}$ as $t\to\infty$, \eqref{proba_collision} follows by choosing $\phi(t) = c \psi(t)^2$ for some  and letting $t$ go to infinity, and then $c\to\infty$. Indeed,
as $t\to\infty$ by the Donsker invariance principle, the first term of the r.h.s. becomes
\[
\pee ( |\caN(0,c)|\leq 1),
\]
where $\caN(0,c)$ denotes a normal random variable with mean zero and variance $c$. While the second term of the r.h.s. vanishes since $\psi(t) = o(\sqrt{t})$, and hence
$\dfrac{|X(\phi(t))|}{\sqrt{t}}\to 0$ as $t\to\infty$. This in turn tends to zero when $c\to\infty$.

At this stage we proved that the probability that the two sets of SIP particles are coupled in the time
lapse $[(1-\delta)t,  t]$ is strictly positive. Then by iterating independently such coupling attempts, one sees that
the probability of eventual successful coupling of the SIP particles is one.
\epr

As a consequence,
\bc\label{cor:success}
Let $\mu$ be a tempered invariant measure for the SIP. Then for all $n\in \N$, $\exists \ \alpha_n\in [0,\infty)$
such that for all $\mathbf x = (x_1,...,x_n) \in \big(\mathbb Z^d\big)^n$,
\[
\widehat \mu (\mathbf x)= \alpha_n
\]
i.e., $\widehat \mu (\mathbf x)$ only depends on $n=|\mathbf x|$ but not on the precise locations $x_1,\ldots,x_n$.
\ec

\begin{proof}
The proof is quite standard, using that bounded harmonic functions are constant when a successful coupling exists, but we give it here for the sake of self-consistency.
Let $\mu$ be a tempered invariant measure, and put, as before, $c_n= \sup_{x_1,\ldots,x_n} \widehat{\mu} (x_1,\ldots, x_n)$.
Consider two sets of SIP-particles $\mathbf X^S(t)=(X_1^S(t),...,X_n^S(t))$ and $ \mathbf Y^S(t)=(Y_1^S(t),...,Y_n^S(t))$. Since there exists a successful coupling for the SIP-particles, for any $\mathbf x=(x_1,...,x_n)$, $\mathbf y=(y_1,...,y_n) \in \big(\mathbb Z^d\big)^n$
$$\widehat {\mathbb P}_{\mathbf x,\mathbf y}^{SIP} ( \tau < \infty ) = 1$$
where $\tau = \inf \big( t \geq 0 : \mathbf X^S(s) = \mathbf Y^S(s) \mbox{ for all } s \geq t \big)$. Now since $\mu \in \mathcal I_t$ by Proposition \ref{prop:harmo},
\begin{align*}
& \widehat \mu (\mathbf x) = \E_{\mathbf x}^{SIP} \widehat \mu (X^S(t)) \\
& \quad =  \widehat {\mathbb E} _{\mathbf x,\mathbf y}^{SIP} \Big( \widehat \mu(\mathbf Y^S (t)) \mathbf 1 \{ \mathbf X^S (t)= \mathbf Y^S (t) \}\Big) + \widehat {\mathbb E} _{\mathbf x, \mathbf y}^{SIP} \Big(\widehat \mu(\mathbf X^S (t)) \mathbf 1 \{ \mathbf X^S (t)\neq \mathbf Y^S (t)\} \Big)
\\
& \quad =
\widehat {\mathbb E} _{\mathbf x,\mathbf y}^{SIP}( \widehat \mu (Y^S(t))
+ \widehat {\mathbb E} _{\mathbf x,\mathbf y}^{SIP}\left(\left(\widehat \mu (X^S(t))- \widehat \mu (Y^S(t))\right)\mathbf 1 \{ \mathbf X^S (t)\neq \mathbf Y^S (t)\}\right)\\
& \quad =
\widehat \mu (\mathbf y) + \widehat {\mathbb E} _{\mathbf x,\mathbf y}^{SIP}\left(\left(\widehat \mu (X^S(t))- \widehat \mu (Y^S(t))\right)\mathbf 1 \{ \mathbf X^S (t)\neq \mathbf Y^S (t)\}\right)
\end{align*}
Hence, using $\widehat{\mu} (x)\leq c_n$ for $x\in (\Zd)^n$, we obtain the estimate
\[
|\widehat \mu (\mathbf x)-\widehat \mu (\mathbf y)|\leq 2c_n \widehat {\mathbb P} _{\mathbf x,\mathbf y}^{SIP}\left(\mathbf X^S (t)\neq \mathbf Y^S (t)\right)
\]
and the result follows by letting $t\to\infty$.
\end{proof}

\section{Ergodic tempered measures}\label{4}
From now on, we can characterize the ergodic tempered measures.
First, in the next lemma we show that this set coincides with the extreme elements of the set of
tempered invariant measures.
\bl\label{hupu}
The set of tempered ergodic invariant measures coincide with the set of extreme points of the tempered invariant measures, i.e.,
\[
(\caI_t)_e= (\caI\cap \caP_t)_e= \caI_e\cap \caP_t
\]
\el
\bpr
Let $\mu\in (\caI\cap \caP_t)$ and $\mu\not\in(\caI\cap \caP_t)_e$, then there exist $0<\lambda<1$ and $\mu_1,\mu_2\in (\caI\cap \caP_t)$ such that
\[
\mu =\lambda \mu_1+ (1-\lambda) \mu_2
\]
therefore $\mu\not\in \caI_e$, and hence $\mu\not\in \caI_e\cap \caP_t$. Conversely, suppose that $\mu\in \caP_t$ is
not in $\caI_e\cap \caP_t$, then $\mu\not\in \caI_e$ and hence
then there exist $0<\lambda<1$ and $\mu_1,\mu_2\in\caI$ such that
\[
\mu =\lambda \mu_1+ (1-\lambda) \mu_2
\]
this equality together with the fact that $\mu\in\caP_t$,  and the positivity of the functions $D(\xi,\cdot)$, imply that $\mu_1,\mu_2\in \caP_t$.
Therefore $\mu\not\in (\caI\cap \caP_t)_e$. Finally if $\mu$ is not  in $\caI_e\cap \caP_t$ and $\mu$ is also not in $\caP_t$ then trivially
$\mu\not\in (\caI\cap \caP_t)_e$.
\epr

We can then characterize the ergodic tempered measures.
\bt\label{temperedinv}
If $\mu\in \caI_e\cap \caP_t$ then $\mu=\nu_\lambda^m$ for some $\lambda\in [0,1)$.
As a consequence
\be\label{clebard}
(\caI_t)_e= \{ \nu_\lambda^m:\lambda\in [0,1)\}.
\ee
\et

\bpr
Remark that by the multivariate version of the Carleman moment condition, within the class $\caP_t$, a measure $\mu$ is uniquely determined by
its $D$-transform $\hat{\mu}$. Let $\mu\in (\caI_t)_e$.
Since $\mu$ is invariant its $D$-transform $\widehat{\mu} (\xi)$ depends only on $|\xi|$ so we put,
with slight abuse of notation $\hat{\mu} (\xi)= \hat{\mu} (n)$.

In order to show that $\mu=\nu_\lambda^m$ for some $\lambda\in [0,1)$, it
suffices now to show that $\hat{\mu}(n)= a^n$ for some $a\geq 0$.
This  in turn follows if we show
that
\be\label{facto}
\hat{\mu}(n+m)= \hat{\mu}(n)\hat{\mu}(m),
\ee
for all $n, m\in\N$.
Denote
$S_t$ the semigroup of the SIP and
denote
$\caS_T=\frac1T\int_0^T S_t dt$.
Fix $\xi, \xi'$ to finite configurations, with $|\xi|=n, |\xi'|=m$.
By ergodicity we have, $\mu$-almost surely
\be\label{waf}
\caS_T D(\xi, \eta)\to \hat{\mu} (\xi)=\hat{\mu}(|\xi|),
\ee
as $T\to\infty$ (where $\caS_T$ works on $\eta$).

Therefore, by dominated convergence,
\be\label{woef}
\int D(\xi', \eta) \caS_T D(\xi, \eta) \mu(d\eta)\to  \hat{\mu} (\xi')\hat{\mu} (\xi).
\ee
On the other hand, by self-duality
\be\label{wif}
\caS_T D(\xi, \eta)=  \sum_{\xi'' \in \Omega_n} \frac{1}{T} \int_0^T p_{t}(\xi, \xi'') dt  D(\xi'', \eta).
\ee
As a finite number of SIP-particles eventually spread out all over the lattice $\Zd$,
for large $T$, the main contribution of the sum over $\xi''$ in the r.h.s. of \eqref{wif}
is from configurations $\xi''$ in which there are no  particles at locations occupied by particles in $\xi'$
(let us denote this property by $\xi'\perp\xi''$).
If this is the case, then $D(\xi'', \eta) D(\xi', \eta)=   D(\xi''+\xi', \eta)$.
Therefore, using that $\mu\in \caP_t$
\be\label{wof}
\sum_{\xi'' \in \Omega_n} \frac{1}{T} \int_0^T p_{t}(\xi, \xi'') dt  D(\xi'', \eta)=
\sum_{\xi'' \in \Omega_n, \xi''\perp\xi'} \frac{1}{T} \int_0^T p_{t}(\xi, \xi'') dt  D(\xi'', \eta) +b_T
\ee
where $b_T\to 0$ as $T\to\infty$.
Therefore,
\begin{eqnarray}\label{wef}
&&a_T + \hat{\mu}(n)\hat{\mu}(m)=\int D(\xi', \eta) \caS_T D(\xi, \eta) \mu(d\eta)
\nonumber\\
&=&\sum_{\xi'' \in \Omega_n, \xi''\perp \xi'} \frac{1}{T} \int_0^T p_{t}(\xi, \xi'') dt  \int D(\xi''+\xi', \eta) \mu(d\eta) + b_T
\nonumber\\
&=& \sum_{\xi'' \in \Omega_n, \xi''\perp \xi'} \left(\frac{1}{T} \int_0^T p_{t}(\xi, \xi'') dt\right) \hat{\mu}(n+m) +b_T
\nonumber\\
&=&
\sum_{\xi'' \in \Omega_n} \frac{1}{T} \left(\int_0^T p_{t}(\xi, \xi'') dt\right) \hat{\mu}(n+m) +b_T + c_T
\nonumber\\
&=&
\hat{\mu}(n+m) +b_T + c_T
\end{eqnarray}
where $a_T\to 0$ by \eqref{woef} and as explained before $b_T, c_T\to 0$.
Letting now $T\to\infty$ gives
\eqref{facto}.

Now combining this with lemma \ref{hupu}, and the fact that all $\nu^m_\lambda$ are elements of $\caI_e\cap \caP_t$, i.e., are ergodic
under the SIP dynamics (see e.g.\ \cite{GRV}) gives the result \eqref{clebard}.
\epr
\br
As a consequence of proposition \eqref{transprop} this result can be transferred to the BEP, showing
that its tempered invariant ergodic measures are product of Gamma distributions, and
to the BMP, showing that its tempered invariant ergodic measures are product of mean zero Gaussians.
\er
Next, we show that all tempered invariant measure satisfy a correlation inequality of the type derived in \cite{GRV}.
\bp
Let $\mu\in\caI_t$, then for all $(x_1,\ldots,x_n)\in \big(\Z^{d}\big)^n$,
\be\label{corin}
\int D\left(\sum_{i=1}^n \delta_{x_i} ,\eta\right)\mu (d\eta)\geq \prod_{i=1}^n\int D\left( \delta_{x_i} ,\eta\right)\mu (d\eta)
\ee
\ep
\bpr
Because every element of $\caI_t$ can be decomposed into extreme elements,
we have $\mu = \int \nu^m_\lambda d\Lambda (\lambda)$, for some probability measure $\Lambda$ on $[0,1)$ and as a consequence, denoting
$\rho(\lambda)= \frac{\lambda}{1-\lambda}$ we have
\begin{align*}
& \int D\Big(\sum\limits_{i=1}^n \delta_{x_i},\eta\Big) d\mu(\eta)  = \int \int D\Big(\sum\limits_{i=1}^n \delta_{x_i},\eta\Big) d\nu^m_\lambda(\eta)  d\Lambda(\lambda)  = \int \rho(\lambda)^n d\Lambda(\lambda) \\
& \qquad \geq \Big( \int \rho(\lambda) \ d\Lambda(\lambda) \Big)^n  = \prod\limits_{i=1}^n \int D( \delta_{x_i},\eta) d\mu(\eta) .
\end{align*}
\epr
\section{Convergence to ergodic product measures}\label{5}

In this section, we give sufficient criteria for a starting measure $\mu$ to converge in the course of time to one of the product measures $\nu^m_\lambda$.
To this purpose, we introduce the following notions of asymptotic independence and homogeneity. We call a function $f: \N^{\Zd}\to\R$ local if it is a finite linear combination of the functions $D(\xi,\cdot)$.

\begin{definition}
We say a measure $\mu$ is asymptotically homogeneous (AH), if there exists $\rho>0$ such that
\begin{equation}\label{homodef}
\lim_{t\to\infty}\sup_x\left| \E_x \int D(\delta_{X(t)}, \eta) \mu(d\eta)-\rho\right|=0.
\end{equation}
\end{definition}
Here $\E_x$ denote expectation w.r.t.\ simple random walk starting at $x$.
Notice that every translation invariant measure with finite moments is trivially AH.
\begin{definition}
We say a measure $\mu$ is asymptotically independent (AI), if for all $n$ and for all choices of local functions $f_1,\ldots,f_n$,
$$\lim\limits_{|y_i-y_j| \rightarrow \infty} \Big( \int \prod\limits_{i=1}^n \tau_{y_i} f_i d\mu - \prod\limits_{i=1}^n \int \tau_{y_i}f_i d\mu \Big) = 0.$$
\end{definition}
Then we have the following result,
\bt \label{3-converge}
Let $\mu$ be tempered, AH and AI, then
\[
\mu S_t\to \nu^m_{\lambda(\rho)}
\]
with $\lambda(\rho)= \frac{\rho}{1+\rho}$.
\et

\bpr
It is equivalent to prove that
\begin{equation}
\lim\limits_{t \rightarrow \infty} \int \mathbb E_\eta D\Big( \sum\limits_{i=1}^n \delta_{x_i}, \eta\Big) d\mu(\eta)  = \int D\Big(\sum\limits_{i=1}^n \delta_{x_i} , \eta\Big) d\nu_{\lambda(\rho)}(\eta).
\end{equation}
Remark the r.h.s is equal to $\rho^n$ thanks to \eqref{def:relation-dualmeasure}. Now, dealing with the l.h.s., rewrite
$$ \int \mathbb E_\eta D\Big( \sum\limits_{i=1}^n \delta_{x_i}, \eta\Big) d\mu(\eta)  = \int \mathbb E_{\mathbf x}^{SIP} D\Big( \sum\limits_{i=1}^n \delta_{X_i^S(t)}, \eta\Big) d\mu(\eta)$$
Using \cite[Lemma 1]{GRV} for ${\mathbf x}\in \big(\mathbb Z^d\big)^n$, it is equal to
$$\mathbb E_{\mathbf x}^{SIP} \int \prod\limits_{i=1}^n D(\delta_{X^S_i(t)} , \eta) d\mu(\eta) + o (t).$$
Indeed, after an arbitrary large time, with probability close to one, the $n$ SIP-particles initially at $\mathbf x$ will have spread out and will be at different locations.
By the Markov property, for any time scale $\psi(t)$ such that $\psi(t)/{t^{1/4}}\to 0 $ as $t\to\infty$
\begin{align*}
 \mathbb E_{\mathbf x}^{SIP} \int \prod\limits_{i=1}^n D(\delta_{x_i} , \eta) d\mu(\eta) & =  \mathbb E_{\mathbf x}^{SIP} \mathbb E_{\mathbf X^S(t-\psi(t))}^{SIP}  \int \prod\limits_{i=1}^n D(\delta_{X_i^S(\psi(t))}, \eta) d\mu)\eta) \\
&   = \mathbb E_{\mathbf x}^{SIP} \mathbb E_{\mathbf X^S(t-\psi(t))}^{IRW} \int \prod\limits_{i=1}^n D(\delta_{X_i^I(\psi(t))}, \eta) d\mu)\eta) + o (t).
\end{align*}
For the last equality, we used the fact that after the large time span $t-\psi(t)$
the SIP-particles are with probability close to one at distance of the order of $\sqrt{t}$ from each other and therefore,
in the remaining time $\psi(t)$ they will not come closer than $\sqrt{t} - t^{1/8}$ to each other, i.e.,
are still far apart and therefore will move as if they are IRW-particles. Therefore, by dominated convergence and since $\mu$ is tempered and AI, the product over $i$ and the integral over $\mu$ can be exchanged at the price of an $o(t)$ term. Once the product is out of the integral, it trivially also comes out of the IRW expectation, and therefore,
\[
\mathbb E_{\mathbf x}^{SIP} \prod\limits_{i=1}^n \mathbb E_{\mathbf X^I_i(t-\psi(t))}^{IRW} \int D(\delta_{X_i^I(\psi(t))},\eta) d\mu(\eta) + o(t).
\]
Use the AH property \eqref{homodef} of $\mu$ to conclude that this expression in turn is equal to
\[
\rho^n + o(t),
\]
we conclude by letting $t\to\infty$.
\epr
\begin{remark}
We can replace the AH and AI assumptions in Theorem \ref{3-converge} by the following assertion:
\be\label{newcondi}
\lim\limits_{t \rightarrow \infty} \sum\limits_y p_t(x,y) \eta(y) = \rho
\ee
where the limit is in $\mu$-probability.

Indeed, assuming \eqref{newcondi} we can write
\begin{align*}
& \mathbb E^{SIP}_{\mathbf x} \mathbb E_{\mathbf X^I(s)}^{IRW} \int \prod\limits_{i=1}^n D(\delta_{X_i^I(t-s)}, \eta) d\mu(\eta) \\
& \qquad = \mathbb E^{SIP}_{\mathbf x} \int \prod\limits_{i=1}^n \sum\limits_{y_i} p_{t-s} (X_i^I(s),y_i) \eta_s(y_i) d\mu(\eta)
\end{align*}
whose r.h.s converges to $\rho^n$ as $t$, then $s$, go to infinity, by \eqref{newcondi} and dominated convergence (because $\mu$ is tempered by assumption).
\end{remark}

%

We conclude with two additional remarks:
\br
The fact that $p$ is nearest neighbor can be replaced without any difficulty by a finite range kernel
(with the same proof, adapting the definition of collision). Presumably it is enough that $p$ is translation invariant and has a finite second moment.
\er
\br
Related to the SIP is the dual KMP process and its generalized so-called ``thermalized'' SIP \cite{CGGR}, where when particles are at nearest neighbor positions, several particles can jump at the same time.
However, if in this thermalized SIP model all the particles are separated (i.e., at distance $>1$), they behave exactly as independent random walkers, and therefore a finite number of them can be successfully coupled just as SIP particles can. This implies that for the thermalized SIP we have the same
set of ergodic tempered measures. The thermalized SIP in turn is the dual process of a mass-redistribution model, called the thermalized BEP in \cite{CGGR}
which generalizes
the KMP (Kipnis Marchioro Presutti) process \cite{KMP}. Hence, the only ergodic tempered measures of this generalized KMP process are also products
of Gamma distributions.
\er

\paragraph*{Acknowledgement} {We thank Pablo A. Ferrari for pointing us to reference \cite{dMP} and
 further valuable comments and motivating discussions.}

\end{document}